\documentclass{article}
\usepackage{amsmath,amssymb,amssymb,bbm,amsthm,graphicx}
\usepackage[round]{natbib}
\usepackage[colorlinks=true,linkcolor=blue,citecolor=blue,urlcolor=blue]{hyperref}

\renewcommand\P{\mathbb P}
\newcommand\E{\mathbb E}

\newcommand\R{\mathbb R}

\newcommand\I{\mathbbm 1}
\newcommand\cF{\mathcal F}
\newcommand\Var{\mathrm{Var}}

\newcommand\Hop[1]{^{H\!op}_{N,M}(#1;\xi)}
\newcommand\Hopa[1]{^{H\!op}_{N,\alpha N}(#1;\xi)}
\newcommand\SK[1]{^{S\!K}_N(#1;J)}
\newcommand\T[1]{^t_{N,M}(#1;\xi,J)}

\renewcommand\[{\begin{equation}}
\renewcommand\]{\end{equation}}

\newcommand\cites[1]{\citeauthor{#1}'s (\citeyear{#1})}

\newtheorem*{maintheorem}{Theorem \ref*{thm:maintheorem}}
\newtheorem*{overlaptheorem}{Theorem \ref*{thm:moments}}
\newtheorem{theorem}{Theorem}

\newtheorem{corollary}[theorem]{Corollary}
\theoremstyle{definition}
\newtheorem{remark}[theorem]{Remark}

\newcommand\proofskip{\vspace{-16pt}}

\begin{document}

{\centering
{\LARGE{The Hopfield Model with Superlinearly Many Patterns}}\\[.2cm]
James~Y.~Zhao\footnote{Department of Mathematics, Stanford University; email \href{mailto:jyzhao@math.stanford.edu}{jyzhao@math.stanford.edu}}\\[.4cm]
\parbox{.77\textwidth}{\textit{Abstract}: We study the Hopfield model where the ratio $\alpha$ of patterns to sites grows large. We prove that the free energy with inverse temperature $\beta$ and external field $B$ behaves like $\beta\sqrt\alpha+\gamma$, where $\gamma=P(\sqrt2\beta,B)$ is the limiting free energy of the Sherrington-Kirkpatrick model with inverse temperature $\sqrt2\beta$ and external field $B$.}

}

\section{Introduction}
The Hopfield model is the system of \textit{configurations} $x\in\{\pm1\}^N$ governed by the \textit{Gibbs measure}
\[G(x)=\frac1{Z(\beta,B)}e^{-\beta H(x)+B\sum_ix_i},\]
where $Z(\beta,B)=\sum_xe^{-\beta H(x)+B\sum_ix_i}$ is a normalising constant to make $G$ a probability measure, $\beta$  and $B$ are constants representing the inverse temperature and external field respectively, and $H$ is the \textit{Hamiltonian}
\[H\Hop x=-\frac1{\sqrt{NM}}\sum_{k=1}^M\big(x\cdot\xi^k\big)^2
=-\frac1{\sqrt N}\sum_{i=1}^N\sum_{j=1}^N\bigg(\frac1{\sqrt M}\sum_{k=1}^M\xi^k_i\xi^k_j\bigg)x_ix_j.\label{eqn:hamiltonian}\]
The special configurations $\xi^k\in\R^N$, $1\le k\le M$, called \textit{patterns}, are themselves random, making the Gibbs measure a random measure; usually they are taken to be vectors of independent Bernoulli variables, but we will also consider more general patterns.

A crucial aspect of any analysis of the Hopfield model is the relationship between the number of patterns $M$ and the number of sites $N$. Most existing results require either $M/N\rightarrow0$ \citep{gentz,boviergayrard96,boviergayrard97} or $M\sim\alpha N$ with $\alpha\le\alpha_0(\beta)$ for some upper bound $\alpha_0$ \citep{boviergayrardpicco,talagrand98,talagrand00}. Our paper focuses on the parameter region $M\sim\alpha N$ for large $\alpha$, or $M/N\rightarrow\infty$.

In order for the problem to make sense for large $\alpha$, our Hamiltonian is normalised by $1/\sqrt{NM}$ rather than the usual $1/2N$; observe that if $1/\sqrt M$ were replaced by $1/\sqrt N$ in (\ref{eqn:hamiltonian}), then the expression in parentheses would diverge as $\alpha\rightarrow\infty$, resulting in a degenerate Gibbs measure. As a result of this normalisation, our inverse temperature $\beta$ differs from the usual inverse temperature by a factor of $\sqrt\alpha/2$.

The key quantity in this parameter region is the \textit{overlap} $S=(x\cdot\xi^1)/\sqrt N$ between the Gibbs configuration $x$ and the first pattern $\xi^1$. \cite{gentz} proved a central limit theorem for the overlap in the case of sublinearly many patterns, while \cite{talagrand00} proved an exponential moment bound for linearly many patterns. In the superlinear case, in particular for $\alpha>4\beta^2$, we are able to replicate Talagrand's moment bound.

Using this bound and the interpolation technique of \cite{guerratoninelli}, we determine the behaviour of the \textit{free energy}
\[F(\beta,B)=\frac1N\log Z(\beta,B)=\frac1N\log\sum_xe^{-\beta H(x)+B\sum_ix_i},\]
in particular relating it to that of the Sherrington-Kirkpatrick (SK) model
\[H\SK x=-\frac1{\sqrt N}\sum_{i=1}^N\sum_{j=1}^NJ_{ij}x_ix_j,\]
where the $J_{ij}$ are independent standard Gaussian variables.

\begin{theorem}
\label{thm:maintheorem}
Suppose the patterns $\xi$ are independent and symmetric with unit variance and bounded eleventh moment. If $\alpha\ge(4+\epsilon)\beta^2$, then for some universal constant $C=C(\epsilon)$,
\[\Big|\E\big[F\Hop{\beta,B}\big]-\E\big[F\SK{\sqrt2\beta,B}\big]-\beta\sqrt\alpha\Big|\le\frac{C\beta^3}{\sqrt\alpha}\Big(1+O\Big(\frac1{\sqrt N}\Big)\Big).\]
\end{theorem}
\vspace{-4pt}

It is known that the SK free energy converges \citep{guerratoninelli} with exponentially high probability to a limit $P(\beta,B)$ \citep{parisi,guerra,talagrand06}. It is also known that the Hopfield free energy concentrates for Bernoulli patterns \citep{talagrand96}, which we generalise in Theorem \ref{thm:concentration}. Thus, we can write down the following behaviour of the limiting free energy, which is strongly supported by the simulation results in Figure \ref{figure}.

\vspace{4pt}
\begin{corollary}
\label{corollary}
For patterns as in Theorem \ref{thm:maintheorem},
\[\lim_{N\rightarrow\infty}F\Hopa{\beta,B}=\beta\sqrt\alpha+P\big(\sqrt2\beta,B\big)+O\big(\beta^3/\sqrt\alpha\big),\]
\end{corollary}
\vspace{-8pt}
where the limit exists in $L^d$ norm for any $d<11/2$.
\vspace{-8pt}
\begin{proof}
Convergence of the mean follows from Theorem \ref{thm:maintheorem} and \cite{guerratoninelli}, while concentration around the mean follows from Theorem \ref{thm:concentration} and Remark \ref{remark}.
\end{proof}

\[\nonumber\]
{\centering
\includegraphics[width=.5\textwidth]{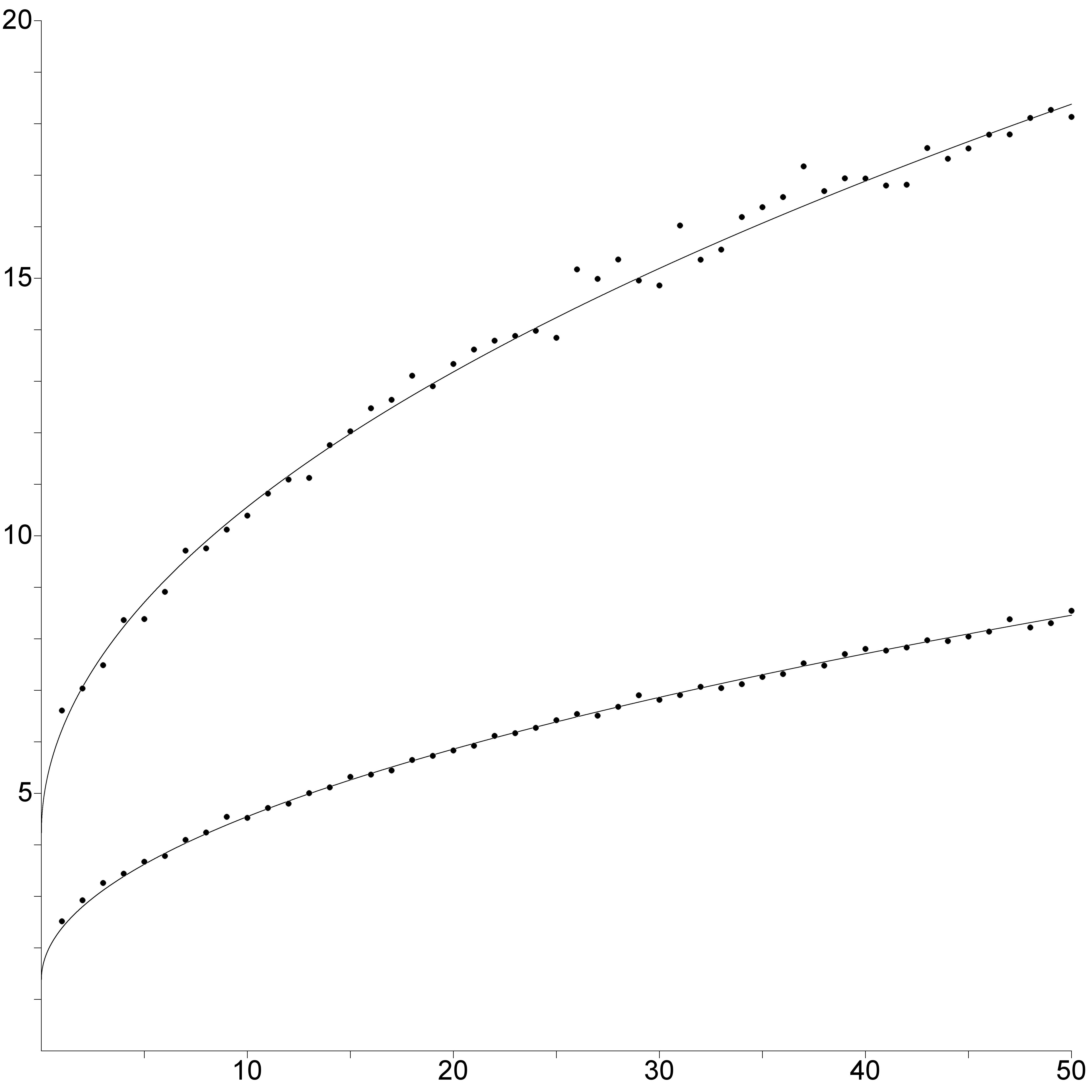}

\parbox{.75\textwidth}{
\textit{Figure 1}. \label{figure}
Plot of free energy (vertical axis) against the parameter $\alpha$ (horizontal axis), for $\beta=1$, $B=0$ (bottom plot) and $\beta=2$, $B=5$ (top plot). The points are realisations of the Hopfield free energy with $N=50$ and $\alpha=1,2,\ldots,50$. The curves are $\beta\sqrt\alpha+P(\sqrt2\beta,B)$, estimated by averaging 100 realisations of $F^{S\!K}_{50}(\sqrt2\beta,B)$.
}

}

\newpage
Heuristically, one can guess that the Hopfield model behaves like the SK model for large $\alpha$, since the interaction between sites $i$ and $j$ is $\frac1{\sqrt M}\sum_k\xi^k_i\xi^k_j$, which converges to Gaussian when $\alpha\rightarrow\infty$ with $N$ fixed. However, this is not easy to prove. Firstly, the Hopfield interactions are not independent, nor can we easily swap the limits $\alpha\rightarrow\infty$ and $N\rightarrow\infty$. Secondly, even ignoring these two issues, convergence of the free energy corresponds to computing the \textit{large deviations} of the Hamiltonian, which are very sensitive; in particular, convergence of the Hamiltonian for typical patterns does not imply corresponding convergence of free energy.

Our strategy is based on the interpolation technique of \cite{guerratoninelli} and the overlap bound of Section \ref{sec:overlap}. First considering Gaussian patterns, we apply \cites{stein} lemma to evaluate the derivative of the interpolating free energy in terms of moments of the overlap, which can be controlled by our overlap bound. Having obtained the result for Gaussian patterns, we extend to Bernoulli and more general patterns following the Taylor expansion method of \cite{carmonahu}; our overlap bound will again be crucial, since the quadratic nature of the Hopfield Hamiltonian as opposed to the linear SK Hamiltonian means that na\"ive bounds are insufficient.

The overlap bound itself is proved using truncation of the distribution function, \cites{hoeffding} inequality and a symmetry trick (\ref{eqn:symmetrytrick}) which allows us to replace the Gibbs measure on configurations with the uniform measure. The symmetry trick (\ref{eqn:symmetrytrick}) may be of interest in itself; it essentially gives a Berry-Esseen type bound for sums of independent random variables when the summands are symmetric, giving a deviation from Gaussian of $O(N^{-p})$ with $p$ as large as desired provided the summands have sufficiently high moments.

Finally, we remark that the simulation results in Figure \ref{figure} fit the curve $\beta\sqrt\alpha+P(\sqrt2\beta,B)$ much closer than the $O(1/\sqrt\alpha)$ error would predict. Computational artefacts and simulation randomness would easily account for the variation, thus the actual error seems to be of much smaller order or perhaps even zero. On the other hand, for $\beta<\frac12$ and $B=0$, the limiting free energy can be written down explicitly \citep{talagrand98}, and the error is precisely $\Theta(1/\sqrt\alpha)$. One might conjecture that a phase transition occurs at $\beta=\frac12$, but we require better control over the error to make further progress.

\textbf{Acknowledgements}. The author wishes to thank Amir Dembo, Persi Diaconis, Jack Kamm and Andrea Montanari for their helpful suggestions which contributed greatly to this paper.

\section{Overlap Bounds}
\label{sec:overlap}
In this section, we show that moments of the \textit{overlap} $S=(x\cdot\xi^1)/\sqrt N$ are bounded uniformly in $N$ and $\alpha$. A similar result appears as Theorem 2.2 in \cite{talagrand00} for linearly many patterns. In the superlinear case, we are able to give a much shorter proof by taking advantage of a symmetry trick. We will require $\alpha>4\beta^2$; observe this does not exclude any interesting cases, since as noted in the introduction, the Gibbs measure becomes degenerate when $\alpha\rightarrow\infty$ but $\beta/\sqrt\alpha$ does not vanish.

\vspace{6pt}
\begin{theorem}[Bernoulli Version]
\label{thm:moments}
Suppose the patterns $\xi$ are Bernoulli. If $\alpha_0>4\beta^2$, then for some $c>0$, $\E\big[e^{cS^2}\big]<\infty$ uniformly in $N$ and $\alpha\ge\alpha_0$.
\end{theorem}
\proofskip
\begin{proof}
Let $E_u=\{u<S^2\le u+1\}$, considered as either a set of configurations or a set of patterns depending on context. Let $H^*=H+S^2/\sqrt\alpha$ be the Hamiltonian with the $\xi^1$ term removed, and let $G^*$ be the corresponding Gibbs measure and $\langle\cdot\rangle^*$ its expectation. Then,
\begin{align}
G\big(S^2>r\big)
&=\sum_{u=r}^\infty G(E_u)
=\sum_{u=r}^\infty\frac{\sum_{x\in E_u}e^{\beta S^2/\sqrt\alpha-\beta H^*\!(x)+B\sum_ix_i}}{\sum_xe^{\beta S^2/\sqrt\alpha-\beta H^*\!(x)+B\sum_ix_i}}\label{eqn:Eu}\\
&\le\sum_{u=r}^\infty e^{\beta(u+1)/\sqrt\alpha}\,\frac{\sum_{x\in E_u}e^{-\beta H^*\!(x)+B\sum_ix_i}}{\sum_xe^{-\beta H^*\!(x)+B\sum_ix_i}}
=\sum_{u=r}^\infty e^{\beta(u+1)/\sqrt\alpha}\,G^*\!(E_u).
\end{align}
Since $G^*$ does not depend on $\xi^1$, we can exchange $\E_{\xi^1}$ and $\langle\cdot\rangle^*$ by Tonelli's theorem to obtain
\[\E\big[G(S^2>r)\big]
\le\sum_{u=r}^\infty e^{\beta(u+1)/\sqrt\alpha}\,\E\big[\big\langle\P_{\xi^1}\!(E_u)\big\rangle^*\big].\label{eqn:tailbound}\]
By \cites{hoeffding} inequality, for any fixed $x$, $\P_{\xi^1}\!(E_u)\le2e^{-u/2}$, hence
\[\E\big[G(S^2>r)\big]
\le2\sum_{u=r}^\infty\exp\Big(\frac{\beta(u+1)}{\sqrt\alpha}-\frac u2\Big)
=Ce^{-ar},\]
where $a=1/2-\beta/\sqrt\alpha>0$ and $C=2e^{\beta/\sqrt\alpha}/(1-e^{-a})$. Since $S^2\ge0$,
\[\E\big[\big\langle e^{cS^2}\big\rangle\big]
=\int_0^\infty\P\Big[G\big(e^{cS^2}>r\big)\Big]dr
=\int_0^\infty\E\Big[G\big(S^2>\tfrac1c\log r\big)\Big]dr
\le1+C\int_1^\infty r^{-a/c}dr.\]
Picking any $0<c<a$ completes the proof.
\end{proof}

Theorem \ref{thm:moments} allows the overlap to be controlled very well and is interesting in its own regard. However, the proof of Theorem \ref{thm:maintheorem} will only require the third moment $\E[\langle|S|^3\rangle]$, which is uniformly bounded for a much larger class of patterns.

\vspace{6pt}
\begin{overlaptheorem}[General Version]
Suppose the patterns $\xi$ are symmetric with variance 1 and bounded $m$th moment for some $m>2$. If $\alpha_0>4\beta^2$ and $0<d<(m-2)(1-d/m)$, then $\E\big[|S|^d\big]<\infty$ uniformly in $N$ and $\alpha\ge\alpha_0$.
\end{overlaptheorem}
\proofskip
\begin{proof}
Pick $\epsilon$ so that $\sqrt{\alpha_0}\ge(2+3\epsilon)\beta$, and let $A=\big\{||\xi^1||_2^2>(1+\epsilon)N\big\}$. Then,
\begin{align}
\E\big[\big\langle|S|^d\big\rangle\big]
&=\E\big[\big\langle|S|^d\big\rangle\I_A\big]+\E\big[\big\langle|S|^d\big\rangle\I_{A^c}\big]\\
&\le\E\big[\big\langle|S|^m\big\rangle\big]^{d/m}\P\big[A\big]^{1-d/m}
+\int_0^\infty\E\Big[G\big(S^2>r^{2/d}\big)\I_{A^c}\Big]dr.\label{eqn:momentsplit}
\end{align}
The bound on the integral is similar to the Bernoulli case; (\ref{eqn:tailbound}) becomes
\[\E\big[G(S^2>r)\I_{A^c}\big]
\le\sum_{u=r}^\infty e^{\beta(u+1)/\sqrt\alpha}\,\E\big[\big\langle\P_{\xi^1}\!(E_u\cap A^c)\big\rangle^*\big].\]
Since replacing $\xi^1_i$ with $x_i\xi^1_i$ changes neither the distribution of $S$ nor the value of $||\xi^1\!||_2$, $\P_{\xi^1}\!(E_u\cap A^c)$ does not depend on $x$, and hence can be replaced by its average over all $x$. Again by \cites{hoeffding} inequality,
\[\P_{\xi^1}\!(E_u\cap A^c)
=\frac1{2^N}\sum_x\P_{\xi^1}(E_u\cap A^c)
=\E\bigg[\frac{|E_u|}{2^N}\I_{A^c}\bigg]
\le\E\big[2e^{-uN/2||\xi^1\!||^2_2}\big]\le 2e^{-u/(2+2\epsilon)}.\label{eqn:symmetrytrick}\]
Since $a=1/(2+2\epsilon)-\beta/\sqrt\alpha>0$ by our choice of $\epsilon$, the same argument shows that the integral in (\ref{eqn:momentsplit}) is bounded. For the remaining term, Cauchy-Schwarz gives $|S|\le||\xi^1||_2$, and by Jensen's inequality, $\E\big[\big\langle|S|^d\big\rangle\big]\le\E\big[||\xi^1||^d_2\big]\le N^{d/2}D$. Thus, we need $\P[A]^{1-d/m}=O(N^{-d/2})$. Since $d<(m-2)(1-d/m)$, it suffices to prove that $\P[A]=O\big(N^{-(m-2)/2}(\log N)^m\big)$.

Let $f(r)=r^m\P\big[|\xi^1_1|>r\big]$. We claim that there exists $b\in\R$ such that for any $r>0$, the interval $[r,br]$ contains a point $s$ with $f(s)<1$. Indeed, if the claim is false, then for any $b$, there is $r_b$ such that $f\ge1$ on $[r_b,br_b]$, hence
\[\frac1m\E\big[|\xi^1_1|^m\big]=\int_0^\infty r^{m-1}\P\big[|\xi^1_1|>r\big]dr=\int_0^\infty\frac{f(r)}rdr\ge\int_{r_b}^{br_b}\frac1rdr=\log b,\]
which is a contradiction since the left side is bounded while the right side is unbounded. Thus, we can pick $r_N$ such that $\sqrt N/\log N\le r_N\le b\sqrt N/\log N$ and $\P\big[|\xi^1_1|>r_N\big]<r_N^{-m}$. Let $B=\big\{\forall i\,|\xi^1_i|\le r_N\big\}$. Then,
\[\P\big[B^c\big]\le N\P\big[|\xi^1_1|>r_N\big]\le Nr_N^{-m}\le N^{1-m/2}(\log N)^m.\]
Finally, by \cites{hoeffding} inequality, $\P\big[B\cap A\big]\le 2e^{-\epsilon N/2r_N^2}\le2e^{-\epsilon(\log N)^2/2b^2}$.
\end{proof}

\section{Main Result}

\begin{maintheorem}[Gaussian Version]
Suppose the Hopfield patterns $\xi$ and SK interactions $J$ are both Gaussian. If $\alpha\ge\alpha_0>\beta^2$, then for some constant $C$ depending only on $\alpha_0/\beta^2$,
\[\Big|\E\big[F\Hop{\beta,B}\big]-\E\big[F\SK{\sqrt2\beta,B}\big]-\beta\sqrt\alpha\Big|\le\frac{C\beta^3}{\sqrt\alpha}.\]
\end{maintheorem}
\proofskip
\begin{proof}
For $0\le t\le1$, define the interpolated Hamiltonian
\begin{align}
H\T x&=\sqrt{1-t}\sqrt2H\SK x+\sqrt tH\Hop x+N\sqrt{\alpha t}\\
&=-\sum_{i,j}\bigg(\frac{\sqrt{1-t}\sqrt2}{\sqrt N}J_{ij}+\frac{\sqrt{t}}{\sqrt{NM}}\sum_k\xi_i^k\xi_j^k\bigg)x_ix_j+N\sqrt{\alpha t}.\label{eqn:Ht}
\end{align}
For the remainder of the proof of Theorem \ref{thm:maintheorem}, Gibbs measure quantites will refer to the interpolated Hamiltonian unless explicitly specified otherwise. Observe that 
\[\E\big[F\Hop{\beta,B}-F\SK{\sqrt2\beta,B}-\beta\sqrt\alpha\big]=\int_0^1\frac d{dt}\E\big[F\T{\beta,B}\big]dt.\]
We can calculate this derivative as
\begin{align}
\frac d{dt}\E[F]
&=\E\bigg[\frac1{NZ}\sum_x\Big({-\beta}\frac{\partial H}{\partial t}\Big)e^{-\beta H(x)+B\sum_ix_i}\bigg]
=-\frac\beta N\E\bigg[\bigg\langle\frac{\partial H}{\partial t}\bigg\rangle\bigg]\\
&=\frac{\beta(N-1)}{\sqrt{2(1-t)N}}\E\big[J_{12}\langle x_1x_2\rangle\big]-\frac{\beta(N-1)M}{2\sqrt{tNM}}\E\big[\xi^1_1\xi^1_2\langle x_1x_2\rangle\big].\label{eqn:deriv}
\end{align}
To obtain (\ref{eqn:deriv}), we differentiate (\ref{eqn:Ht}) and split the sum over and $i$ and $j$ into two sums where $i\ne j$ and $i=j$ respectively. The $i\ne j$ terms give (\ref{eqn:deriv}) upon taking $i=1$, $j=2$ and $k=1$ without loss of generality. For the $i=j$ terms, the SK term is a multiple of $\E\big[J_{11}\langle x_1^2\rangle\big]=0$, while the Hopfield term exactly cancels the derivative of $N\sqrt{\alpha t}$.

Let $X=x_1x_2$ and $\Var(X)=\langle X^2\rangle-\langle X\rangle^2$. By \cites{stein} lemma,
\[\E\big[J_{12}\langle x_1x_2\rangle\big]
=\E\bigg[\frac\partial{\partial J_{12}}\langle x_1x_2\rangle\bigg]
=\frac{\beta\sqrt{2(1-t)}}{\sqrt N}\E\big[\Var(X)\big].\label{eqn:SKterm}\]
For any function $g=g(x;\xi,J)$,
\[\frac\partial{\partial\xi^1_i}\big\langle g\big\rangle=\bigg\langle\frac{\partial g}{\partial\xi^k_i}\bigg\rangle+\frac{2\beta\sqrt t}{\sqrt M}\Big(\big\langle gx_iS\big\rangle-\big\langle g\big\rangle\big\langle x_iS\big\rangle\Big).\label{eqn:diffrule}\]
Let $U=x_1S$ and $V=x_2S$. By \cites{stein} lemma and (\ref{eqn:diffrule}),
\begin{align}
\E\big[\xi^1_1\xi^1_2\langle x_1x_2\rangle\big]\label{eqn:hopfieldterm}
&=\frac{2\beta\sqrt t}{\sqrt M}\E\bigg[\frac\partial{\partial\xi^1_1}\Big(\langle U\rangle-\langle X\rangle\langle V\rangle\Big)\bigg]\\
&=\frac{2\beta\sqrt t}{\sqrt{NM}}\E\big[\Var(X)\big]+\frac{4\beta^2t}{M}\E\Big[\langle S^2\rangle-\langle U\rangle^2-\langle V\rangle^2-\langle X\rangle\langle UV\rangle+2\langle X\rangle\langle U\rangle\langle V\rangle\Big].\nonumber
\end{align}
The first term exactly cancels (\ref{eqn:SKterm}), while the second term is $O(\beta^2/M)$ by H\"older's inequality and Theorem \ref{thm:moments}, which gives the $C\beta^3/\sqrt\alpha$ term.
\end{proof}

To extend the result from Gaussian patterns to general patterns, we follow the Taylor expansion strategy of \cite{carmonahu}. However, we note that our overlap bound is crucial for this to work; naive bounds are insufficient as the Hopfield Hamiltonian is quadratic in the patterns, while the SK Hamiltonian is linear in the interactions.

\vspace{6pt}
\begin{maintheorem}[General Version]
Suppose the patterns $\xi$ are symmetric with unit variance and bounded eleventh moment. If $\alpha\ge\alpha_0>\beta^2$, then for some $C$ depending only on $\alpha_0/\beta^2$,
\[\Big|\E\big[F\Hop{\beta,B}\big]-\E\big[F\SK{\sqrt2\beta,B}\big]-\beta\sqrt\alpha\Big|\le\frac{C\beta^3}{\sqrt\alpha}+O\Big(\frac{\beta^3}{\sqrt M}\Big).\]
\end{maintheorem}
\proofskip
\begin{proof}
The Gaussian assumption is only used in (\ref{eqn:hopfieldterm}), so it suffices to prove
\[\bigg|\E\big[\xi^1_1\xi^1_2\langle x_1x_2\rangle\big]
-\E\bigg[\frac\partial{\partial\xi^1_1}\frac\partial{\partial\xi^1_2}\langle x_1x_2\rangle\bigg]\bigg|=O\Big(\frac{\beta^2}{M\sqrt N}\Big).\label{eqn:steinsmethod}\]
For readers unfamiliar with \cite{carmonahu}, we note that while (\ref{eqn:steinsmethod}) cosmetically resembles Stein's method, we are not proving that $\xi^1_1$ and $\xi^1_2$ become Gaussian (by definition they do not), but rather, that the dependence of $\langle x_1x_2\rangle$ on $\xi^1_1$ and $\xi^1_2$ vanishes, in the sense that \cites{stein} lemma holds for non-Gaussian variables with some small error.

Define the probabilistic function $f:\R^2\rightarrow\R$ by $f=\langle x_1x_2\rangle$ considered as a function of $y=\xi^1_1$ and $z=\xi^1_2$, with the other $\xi^k_i$ regarded as random parameters. Since $f$ is infinitely differentiable, we can define its second degree Taylor expansion $T$ around 0, that is,
\[T(y,z)=f(0,0)+yf_y(0,0)+zf_z(0,0)+\frac12y^2f_{yy}(0,0)+\frac12z^2f_{zz}(0,0)+yzf_{yz}(0,0).\]
Since the evaluated derivatives depend only on the $\xi^k_i$ other than $y$ or $z$, which are independent with mean 0 and variance 1, an easy calculation yields $\E\big[yzT(y,z)\big]=\E\big[f_{yz}(0,0)\big]$, hence
\[\E\Big[yzf(y,z)-f_{yz}(y,z)\Big]=\E\Big[yz(f-T)(y,z)-\big(f_{yz}(y,z)-f_{yz}(0,0)\big)\Big].\]
By Taylor's theorem,
\begin{align}
(f-T)(y,z)&=\int_0^1\Big(\tfrac12y^3f_{yyy}(sy,sz)+\tfrac32y^2zf_{yyz}(sy,sz)\nonumber\\
&\hspace{1cm}{}+\tfrac32yz^2f_{yzz}(sy,sz)+\tfrac12z^3f_{zzz}(sy,sz)\Big)(1-s)^2ds;\\
f_{yz}(y,z)-f_{yz}(0,0)&=\int_0^1\Big(yf_{yyz}(sy,sz)+zf_{yzz}(sy,sz)\Big)ds.
\end{align}
By symmetry of $f(y,z)=f(z,y)$,
\begin{align}
\E\big[yzf(y,z)-f_{yz}(y,z)\big]
&=\int_0^1(1-s^2)\E\big[y^4zf_{yyy}(sy,sz)\big]\nonumber\\
&\hspace{1.3cm}{}+\E\big[\big(3(1-s)^2y^3z^2+2y\big)f_{yyz}(sy,sz)\big]ds.\label{eqn:taylortheorem}
\end{align}
Let $X=x_1x_2$ and $U=x_1S$. By repeated application of (\ref{eqn:diffrule}),
\begin{align}
f_y&=\bigg(\frac{2\beta\sqrt t}{\sqrt M}\bigg)\phantom{{}^2}\Big(\langle XU\rangle-\langle X\rangle\langle U\rangle\Big);\label{eqn:firstderiv}\\
f_{yy}&=\bigg(\frac{2\beta\sqrt t}{\sqrt M}\bigg)^2\Big(\langle XU^2\rangle-2\langle XU\rangle\langle U\rangle-\langle X\rangle\langle U^2\rangle+2\langle X\rangle\langle U\rangle^2\Big);\label{eqn:secondderiv}\\
f_{yyy}&=\bigg(\frac{2\beta\sqrt t}{\sqrt M}\bigg)^3\bigg(
\langle XU^3\rangle-3\langle XU^2\rangle\langle U\rangle-3\langle XU\rangle\langle U^2\rangle+6\langle XU\rangle\langle U\rangle^2\nonumber\\
&\hspace{4cm}{}-\langle X\rangle\langle U^3\rangle+6\langle X\rangle\langle U^2\rangle\langle U\rangle-6\langle X\rangle\langle U\rangle^3
\bigg).\label{eqn:thirdderiv}
\end{align}
By H\"older's inequality, $\big|f_{yyy}\big|\le8\beta^3M^{-3/2}26\big\langle|S|^3\big\rangle$. Let $\langle\cdot\rangle^s$ be expectation with respect to the Gibbs measure with $\xi^1_1$ and $\xi^1_2$ replaced by $s\xi^1_1$ and $s\xi^1_2$ respectively. Then,
\[\Big|\E\big[y^4zf_{yyy}(sy,sz)\big]\Big|
\le\frac{8\beta^3}{M^{3/2}}26\,\E\bigg[|y|^4|z|\Big\langle\big|S-(1-s)(x_1y+x_2z)\big|^3\Big\rangle\!\raisebox{5pt}{${}^s$}\bigg].\]
Thus, we need to bound $\E\big[\big\langle|y|^{4+a}|z|^{1+b}|S|^c\big\rangle\!\raisebox{2pt}{${}^s$}\big]$, where $a+b+c=3$. By H\"older's inequality, this is bounded by $\E\big[|y|^{(4+a)d/(d-c)}|z|^{(1+b)d/(d-c)}\big]{}^{1-c/d}\,\E\big[\big\langle|S|^d\big\rangle\big]{}^{c/d}$, where $(4+a)d/(d-c)\le m$ and $c<d<(m-2)(1-d/m)$. By Theorem \ref{thm:moments}, which applies to $\langle\cdot\rangle^s\vphantom{\big|}$ since scaling finitely many patterns does not affect the proof, this yields $\E\big[y^4zf_{yyy}(sy,sz)\big]=O\big(\beta^3/M^{3/2}\big)$; some algebra shows we need $m>6+\sqrt{22}\approx10.7$.

Bounding the second expectation in (\ref{eqn:taylortheorem}) is similar. Using (\ref{eqn:diffrule}) to differentiate (\ref{eqn:secondderiv}) with respect to $z$ and applying H\"older's inequality gives
\[\big|f_{yyz}\big|\le\frac{4\beta^2}{M\sqrt N}8\big\langle|S|\big\rangle+\frac{8\beta^3}{M^{3/2}}26\big\langle|S|^3\big\rangle.\]
Scaling $\xi^1_1$ and $\xi^1_2$ by $s$, multiplying by $3(1-s)^2y^3z^2+2y$ and taking expectations, we obtain the bound $O(\beta^2/M\sqrt N+\beta^3/M^{3/2})$, which is $O(\beta^2/M\sqrt N)$ since $\beta/\sqrt\alpha<\frac12$.
\end{proof}

Finally, we prove the concentration that allows us to conclude Corollary \ref{corollary} from Theorem \ref{thm:maintheorem}. We follow the martingale approach of \cite{pasturshcherbina} and \cite{carmonahu}, again relying on our overlap bound to complete the proof.

\vspace{6pt}
\begin{theorem}
\label{thm:concentration}
Suppose $\alpha\ge(4+\epsilon)\beta^2$ for some $\epsilon>0$, and the patterns $\xi$ have bounded $m$th moment. If $4\le2d<(m-2)(1-2d/m)$, then for some $C$ depending only on $\epsilon$ and $\E\big[|\xi^1_1|^m\big]$,
\[\E\bigg[\Big|F\Hopa{\beta,B}-\E\big[F\Hopa{\beta,B}\big]\Big|^d\bigg]\le\frac{C\beta^d}{N^{d/2}}.\]
\end{theorem}
\proofskip
\begin{proof}
Let $D_\ell=N\big(\E[\log Z|\cF_\ell]-\E[\log Z|\cF_{\ell-1}]\big)$, where $\cF_\ell$ is the $\sigma$-algebra generated by $\xi^1,\ldots,\xi^\ell$. The $D_\ell$ are martingale differences with $\sum_\ell D_\ell=N(F-\E[F])$, so by \cites{burkholder} martingale inequality and convexity of $x\mapsto x^{d/2}$, for some universal constant $C_d$,
\[\E\big[|F-\E[F]|^d\big]=\frac1{N^d}\E\bigg[\Big|\sum_\ell D_\ell\Big|^d\bigg]
\le\frac{C_d}{N^d}\E\bigg[\Big(\sum_\ell D_\ell^2\Big)^{d/2}\bigg]
\le\frac{C_dM^{d/2}}{N^d}\frac1M\sum_\ell\E\big[|D_\ell|^d\big].\]
Let $\beta'=\beta/\sqrt{NM}$. Observe that for any fixed $1\le\ell\le M$, we can write
\begin{align}
\log Z
&=\log\sum_xe^{\beta'\sum_{k\ne\ell}(x\cdot\xi^k)^2+B\sum_ix_i}-\log\frac{\sum_xe^{\beta'\sum_{k\ne\ell}(x\cdot\xi^k)^2+B\sum_ix_i}}{\sum_xe^{\beta'\sum_k(x\cdot\xi^k)^2+B\sum_ix_i}}\\
&=\log\sum_xe^{\beta'\sum_{k\ne\ell}(x\cdot\xi^k)^2+B\sum_ix_i}-\log\Big\langle e^{-\beta'(x\cdot\xi^\ell)^2}\Big\rangle.
\end{align}
Since the first term above is independent of $\xi^\ell$, it follows that
\[D_\ell=-\E\Big[\log\Big\langle e^{-\beta'(x\cdot\xi^\ell)^2}\Big\rangle\Big|\cF_\ell\Big]+\E\Big[\log\Big\langle e^{-\beta'(x\cdot\xi^\ell)^2}\Big\rangle\Big|\cF_{\ell-1}\Big].\]
By H\"older's inequality,
\[\E\big[|D_\ell|^d\big]
\le2^d\,\E\bigg[\Big|\log\Big\langle e^{-\beta'(x\cdot\xi^\ell)^2}\Big\rangle\Big|^d\bigg]
=2^d\,\E\bigg[\Big|\log\Big\langle e^{-\beta S^2/\sqrt\alpha}\Big\rangle\Big|^d\bigg].\]
The function $\phi(z)=|\log z|^d$ has second derivative $\phi''(z)=\big(d(d-1)\,|\!\log z|^{d-2}-d\,|\!\log z|^{d-1}\big)/z^2$, so $\phi$ is convex for $0<z\le e^{d-1}$. Since $0<e^{-\beta'S^2}\le1<e^{d-1}$, by Jensen's inequality,
\[\E\big[|D_\ell|^d\big]
\le2^d\,\E\bigg[\bigg\langle\Big|\log\big(e^{-\beta S^2/\sqrt\alpha}\big)\Big|^d\bigg\rangle\bigg]
=2^d\,\E\bigg[\bigg\langle\Big(\frac{\beta S^2}{\sqrt\alpha}\Big)^d\bigg\rangle\bigg].\]
This is $O(\beta^d/\alpha^{d/2})$ by Theorem \ref{thm:moments}.
\end{proof}

\begin{remark}
\label{remark}
In the context of our main result, we assume the eleventh moment is bounded, so $m=11$, in which case $4\le2d<(m-2)(1-2d/m)$ if and only if $2\le d<11/2$. In particular, given the assumptions of Theorem \ref{thm:maintheorem}, Theorem \ref{thm:concentration} shows concentration around the mean in $L^d$ norm for any $d<11/2$, as claimed in Corollary \ref{corollary}.
\end{remark}

\newpage
\newcommand\MR[1]{\href{http://www.ams.org/mathscinet-getitem?mr=#1}{MR#1}.}
\setlength\bibhang{1cm}

\end{document}